# Logistics Hub Capacity Deployment in Hyperconnected Transportation Network Under Uncertainty

### Abstract ID: 3426

Xiaoyue Liu, Jingze Li, Benoit Montreuil

*Physical Internet Center, H. Milton Stewart School of Industrial & Systems Engineering,*
*Georgia Institute of Technology, Atlanta, GA 30332*

## Abstract

Modern logistics systems worldwide are facing unprecedented challenges due to the explosive growth of e-commerce, driving the need for resilient systems to tackle problems such as vulnerable supplies, volatile demands, and fragile transportation networks. Motivated by the innovative concept of the Physical Internet, this paper focuses on resilient capacity deployment of open-access logistics hubs in hyperconnected transportation under demand uncertainty and geographical disruptions. We propose a two-stage stochastic optimization model, aiming to smartly deploy the hub capacity to achieve delivery timeliness, high consolidation and network resilience while minimizing hub set-up budget and truck fleet cost. Four optimal hub network configurations are derived by applying scenarios at four stress testing levels into the optimization model, including deterministic demands without hub disruptions, deterministic demands with hub disruptions, stochastic demands without hub disruptions as well as stochastic demands with hub disruptions. To test the performances of different optimal networks, a simulation-based study is then performed over an automotive delivery-to-dealer network and dataset in the Southeast US region. Our results demonstrate the key impacts of various uncertainties on hub capacity deployment in terms of capacity configuration distribution, network resilience, delivery timeliness, and cost-effectiveness. Overall, this study provides a reliable network capacity deployment approach with persistent and sustainable economic and social performances in hyperconnected networks, and the results validate the relationship between capacity deployment and network resilience under different types of uncertainties.

## 1. Introduction

The rapid growth of e-commerce brings challenges to traditional transportation networks, as it boosts the demands for logistics services and introduces complexities in policy and society management. Those issues aggravate temporal and spatial logistic uncertainties as shown in Fig.1. Uncertainties such as unpredictable demand fluctuations and geographical disruptions lead to shipment delays, traffic congestion, economic losses, and resource wastage, further exacerbating the vulnerability of transportation networks. As a result, effective adaptation of transportation networks to cope with uncertainty has risen in importance. Transportation hub networks are strategically designed systems comprised of logistics hubs that serve as focal points for the movement of freight and/or passengers. Their high freight consolidation and flexible origin-destination (OD) connectivity can improve delivery efficiency and reducing transport costs beyond traditional end-to-end transportation [1]. We here focus on enhancing how they cope with uncertainty.

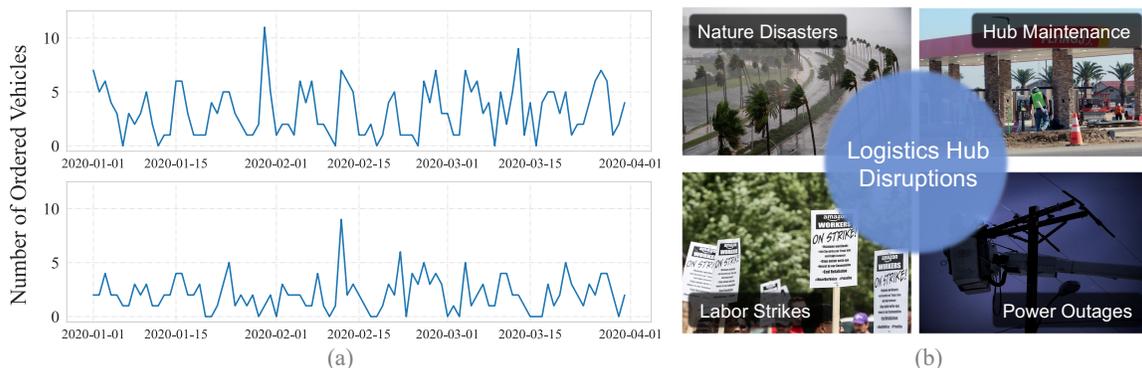

Figure 1: Example (a) shows fluctuation of demands and (b) displays some primary causes of hub disruptions

Given an existing network and the potential locations of open-access logistics hubs, the scope of this paper is to strategically allocate the capacity of each hub to enhance the resilience of transportation networks and better serve the logistic needs under uncertainty. This paper primarily delves into the throughput capacity of logistics hubs, with a particular emphasis on crucial factors such as the quantity of open loading docks, the availability of outbound drivers, the equipment in use and so forth. Unlike physical capacity (determined by hub size, layout, and infrastructure), this capacity is flexible and can be adjusted through planning. We further concentrate on the hyperconnected relay transportation network, where long-haul freight shipments are to be shipped via consecutive short-haul route segments through open-access logistics hubs in a relay way. In the paper, we present a two-stage stochastic optimization model to address spatial-temporal uncertainty, which can be efficiently solved using the Sample Average Approximation (SAA) method. To evaluate the performance of our model and analyze the effects of uncertainty types on logistics hub capacity allocation, we test our optimal hub capacity deployment solutions through a large-scale automotive vehicle delivery simulation study under different types of uncertainty scenarios.

The rest of this paper is organized as follows. Section 2 illustratively describes the problem setting and states research assumptions. In Section 3, we summarize literature reviews of the most relevant papers. Section 4 formulates the capacity deployment optimization model based on two-stage stochastic programming. Section 5 presents a simulation-based study and analyzes experimental results. Section 6 concludes our findings and discusses future research avenues.

## 2. Problem Description

This paper focuses on the resilient capacity deployment problem of open-access logistics hubs for long-haul freight, which require to be transferred by short-haul drivers through multiple hubs in hyperconnected transportation networks. In our general problem setting, given the demand orders, freights are to be shipped from their sources (i.e., origins) to dealer clusters (i.e., destinations) via hubs. In the specific context of vehicle delivery, freight is a car, a utilitarian vehicle, etc.; destination is a dealer or a compact cluster of dealers; source is a plant, port, or railhead used by multiple Original Equipment Manufacturers (OEMs); and open-access hubs in hyperconnected networks enable demand flow through functions such as consolidation, sorting and transfer, providing enhanced hub accessibility and connectivity.

In the simulation study, we use an automotive finished vehicle network in the southeastern USA. Fig. 2 (a) depicts the network nodes consisting of 24 open-access logistics hubs, 13 sources and 50 dealer clusters. Fig. 2 (b) depicts the hyperconnected transportation network with its edges designed to keep the travel time along each edge to less than 5.5 hours, allowing for a short-haul driver round trip to be completed within 11 hours, respecting driving limits imposed by US driving regulations. Demand flow is across the southeast region as shown in Fig. 2 (c) and fluctuates daily according to customer needs. We consider two types of uncertainties: demand uncertainty and disruption uncertainty. Disruptions randomly occur at hubs and cause significant delays in travel time on the arcs connecting to the affected hubs. Given these settings, we aim to determine hub capacity in this hyperconnected transportation network to serve logistics demands while minimizing hub costs and fleet costs under various uncertainty scenarios.

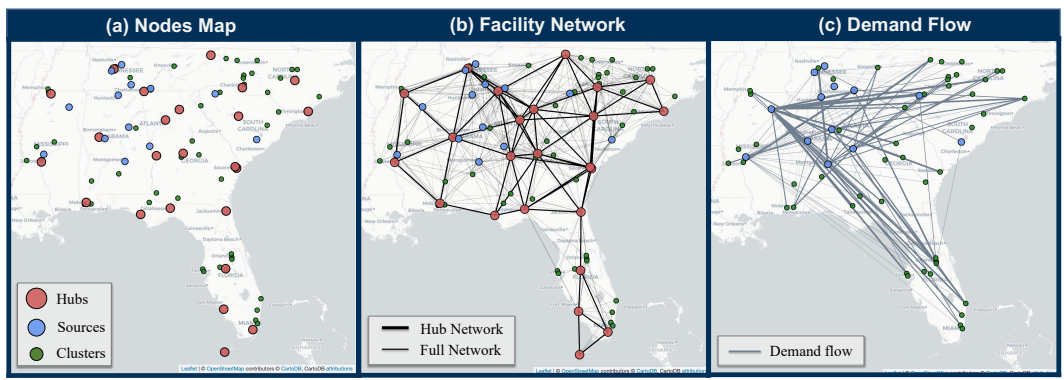

Figure 2: Nodal map (a), facility network (b) and demand flow (c) for multiple automotive OEMs in Southeast USA

## 3. Related Research

Much work has been done on exploring network configurations for logistic and transportation systems. The earliest and simplest configuration is the direct node-to-node linkages from origins to destinations. However, its discrete commodity routing, handling, and sorting leads to massive network operating costs. To mitigate this problem, the hub-and-spoke configuration has been introduced to leverage scale economies by consolidating freight flows [2]. The hub-and-spoke network is widely adopted in transportation, and numerous studies have been developed based on this configuration. However, its centralized structure causes inefficient detours, extended delivery time, and potential network crashes from hub disruptions. None of these issues could be solved intrinsically until the hyperconnected network configuration was initially proposed as a key Physical Internet concept [3]. In contrast to hub-and-spoke networks, hyperconnected networks are multi-plane interconnected meshed hub networks that can smoothly support dynamic freight routing and hub-based consolidation and strengthen connectivity between logistics hubs [4]. By taking advantage of modularity, containerization, and hyperconnectivity, these networks can offer faster and more reliable deliveries from origins to destinations, thereby improving efficiency, service capability, resilience, and sustainability.

Various studies have been conducted to enable and leverage resilient network design. [5] develops a bi-level, three-stage stochastic optimization model to tackle network design for disaster management including pre-event mitigation, preparedness, and post-event response. In [6], a trilevel optimization model is proposed for the resilient network design problem to minimize the total travel time under uncertain disruptions. When considering the resilience of the Physical Internet, its hyperconnectivity tends to provide more efficient and resilient networks compared to traditional networks [7]. Thus, studying the full potential of the Physical Internet to enhance network resilience is a valuable research area.

In the field of capacity deployment, numerous papers have tackled the capacity facility location problem (CFLP) for hub networks. [8,9] formulate the hub location problem under capacity constraints and optimize the size of hubs by determining the discrete capacity levels for each hub. In [10], a dynamic deployment optimization model is proposed to address the modular smart lockers relocation problem in hyperconnected logistics systems, taking into account uncertain events. The model comprises a two-stage stochastic optimization approach, and to efficiently solve large-scale instances of this problem, Bender's decomposition, combined with an acceleration method, is developed in [10].

In this paper, we aim to fill the gaps of the resilient capacity hub network design in hyperconnected transportation by optimizing the capacity deployment of logistics hub to enhance network resilience. Compared to the present state of research, we make the following contributions. First, our study explores the potential of hyperconnected network resilience through a concrete methodology and simulation analysis. Second, we validate that the deployment of logistics hub capacities has a significant impact on the transportation network resilience, opening avenues for future research on resilient capacity deployment. Third, the hub networks resulting from our work can become inputs for other studies based on the Physical Internet, such as operating system design [11] and simulation modeling [12].

Table 1: Notations of sets, parameters, and decision variables in the two-stage stochastic model

| Sets | | Parameters | |
|---|---|---|---|
| $O$ | set of origins, indexed by $o$ | $q_{od}(\omega)$ | demands from origin $o$ to destination $d$ |
| $D$ | set of destinations, indexed by $d$ | $t(a,\omega)$ | travel time on arc $a$ in scenario $\omega$ |
| $H$ | set of logistics hub, indexed by $h$ | $s_h$ | hub operational cost per capacity $h$ |
| $V$ | set of nonempty demand $OD$ pairs | $f_h$ | fixed cost of setting up a new hub $h$ |
| $A$ | set of arcs, indexed by $a$ | $v_a$ | fleet cost per hour on arc $a$ |
| $W$ | set of scenarios, indexed by $\omega$ | $m$ | truckload |
| $\delta^+(i)$ | set of arcs come out from node $i$ | $b$ | maximal capacity based on the budget limit |
| $\delta^-(i)$ | set of arcs come into node $i$ | | |
| **First-stage Decision variables** | | **Second-stage Decision variables** | |
| $X_h$ | binary variable, whether to activate hub $h$ or not | $F_o^d(a,\omega)$ | amount of freight from origin $o$ to destination $d$ on arc $a$ in scenario $\omega$ |
| $C_h$ | capacity to deploy at logistics hub $h$ | $Tr(a,\omega)$ | number of trucks on arc $a$ in scenario $\omega$ |

## 4. Methodology

We propose a two-stage stochastic optimization model that aims to minimize the hub capacity costs and the fleet costs over uncertainty scenarios. In the framework of the two-stage stochastic programming, the first stage is to make hub capacity decisions before freight routing schedules are established, while the second stage is to evaluate the total fleet cost due to a realized network configuration, considering the impact of network disruptions and delay penalties on the chosen arcs. The model notations are displayed in Table 1 and its objective function and constraints are as follows:

$$\text{Minimize} \quad \sum_{h \in H} (f_h \cdot X_h + s_h \cdot C_h) + E_\omega \left[ \sum_{a \in A} v_a \cdot t(a,\omega) \cdot Tr(a,\omega) \right] \quad (1)$$

$$\text{subject to} \quad \sum_{a \in \delta^+(o)} F_o^d(a,\omega) = q_{od}(\omega) \qquad \forall (o,d) \in V, \forall \omega \in W \quad (2)$$

$$\sum_{a \in \delta^-(d)} F_o^d(a, \omega) = q_{od}(\omega) \qquad \forall (o,d) \in V, \forall \omega \in W \quad (3)$$

$$\sum_{a \in \delta^+(h)} F_o^d(a, \omega) = \sum_{a \in \delta^-(h)} F_o^d(a, \omega) \qquad \forall h \in H, \forall (o,d) \in V, \forall \omega \in W \quad (4)$$

$$\sum_{(o,d) \in V} F_o^d(a, \omega) \leq m \cdot Tr(a, \omega) \qquad \forall a \in A, \forall \omega \in W \quad (5)$$

$$\sum_{(o,d) \in V} \sum_{a \in \delta^-(h)} F_o^d(a, \omega) \leq C_h \qquad \forall h \in H, \forall \omega \in W \quad (6)$$

$$C_h \leq b \cdot X_h \qquad \forall h \in H \quad (7)$$

$$C_h \in \mathbb{N}, F_o^d(a, \omega) \in \mathbb{N}, Tr(a, \omega) \in \mathbb{N}, X_h \in \{0,1\} \; \forall h \in H, \forall a \in A, \forall (o,d) \in V, \forall \omega \in W \quad (8)$$

The first term in the objective function (1) represents the first-stage hub costs, including fixed hub setup cost and daily operational cost. The second term is the expected value of the second-stage fleet costs, where $\sum_{a \in A} v_a \cdot t(a, \omega) \cdot Tr(a, \omega)$ refers to the total fleet costs under scenario $\omega$. In the formulation, $v_a$, $s_h$ and $f_h$ are cost parameters based on realistic factors such as the cost of land, regional average salary, and the level of local economic development. Flow conservation constraints (2) − (4) ensure that each freight demand is routed from its source to its destination. Constraints (5) provide the minimal number of trucks needed on each arc given the maximal loading capacity of the truck, while constraints (6) prevent capacity restriction violation for each logistics hub. Constraints (7) guarantee that logistics hubs are open when they have active capacity. Constraints (8) specifies the domains of the model variables.

The objective function (1) is computationally intractable due to integration over the continuous distribution of uncertainty scenarios. However, it can be approximated by the formula (9) using the SAA method.

$$\text{Minimize} \quad \sum_{h \in H} ( f_h \cdot X_h + s_h \cdot C_h) + \frac{1}{|W|} \sum_{\omega \in W} ( \sum_{a \in A} v_a \cdot t(a, \omega) \cdot Tr(a, \omega) ) \quad (9)$$

The true distributions of daily demands and random disruption are estimated by the empirical distribution of the $|W|$ independent and identically distributed samples. As $|W|$ approaches infinity, the optimal value and solutions of the model using objective (9) converge to the optimal value and solutions of the original one under mild assumptions [13].

We effectively handle uncertainties by utilizing the two-stage stochastic optimization approach. Regarding demand uncertainty, the demands between each OD pair changes in each scenario $\omega$, as reflected in $q_{od}(\omega)$. The model finds a feasible hub capacity solution to meets demands in all scenarios by satisfying constraints (2) to (4) (network flow constraints). In terms of disruption uncertainty, a hub causes local delays during disruption, leading to the travel time of arcs connecting to this hub increases, as reflected in $t(a, \omega)$. Such increased travel time will result in a rise in fleet costs, which represents a delay penalty for the affected arc. Then the model avoids using those disrupted arcs as we minimize the total fleet costs in objective function (9).

## 5. Simulation-Based Stress Testing Investigation

In this section, we apply our two-stage stochastic optimization model to the Southeast US finished vehicle delivery network presented in Section 2. The synthetic dataset includes around 35,000 vehicle orders over a 3-month planning horizon. In our Southeast US dataset, we observe that utilizing this hyperconnected transportation network can enable all shipments to be delivered within a single day, without carrying over to the next day's transportation resources. Our model exploits this observation to reduce model complexity by avoiding the introduction of a temporal dimension. To generate uncertainty scenarios from the dataset, we estimate the distribution of daily demands on each OD pair using the Kernel Density Estimation (KDE) approach and then take 50 samples from the estimated distribution as its demand scenarios. As for the disruption scenarios, we randomly generate 50 scenarios by assuming that open-access logistics hubs are independent of each other and the disruption rate for each hub $h$ is given by a parameter $p_h$.

To examine the separated effects of demand uncertainty and random geographical disruptions on capacity deployment, we stress test our model at four stress testing levels. The first level does not consider any uncertain conditions, the second level assumes there is no disruption in the network but account for uncertain demands, the third level only accounts for random disruptions given certain demand information, and the fourth level considers both types of uncertainties. Then we achieve these four stress testing levels by inputting four types of scenarios into our model as follows. In the first level, we solve our model using the same 50 scenarios with a 70-percentage daily demands and no disruption in any of the scenarios. Being that all scenarios are identical; this level is equivalent to a deterministic optimization model. For the second level, we use our 50 generated demand scenarios with no disruption in scenarios. In the third level, where we only consider random disruptions, we use our 50 generated disruption scenarios with the 70-percentage daily demands. In the final level, which includes both types of uncertainties, we use the 50 demand scenarios combined with the generated disruption scenarios as our input scenario set.

Solving our two-stage stochastic optimization model at these four stress testing levels yields four optimal networks: the baseline deterministic network (BDN), the stochastic demand network (SDN), the stochastic disruption network (SDiN), and the integrated stochastic network (ISN) as presented in Fig. 3. Combining Table 2 and Fig. 3, we can see that the BDN and SDN tend to allocate large capacity to a small group of logistics hubs located in central areas. It reasonably reflects that when there is no disruption in the network, utilizing few central hubs leads to shorter shipment routes and improved consolidation efficiency, resulting in a decrease in both hub costs and fleet costs. In addition, demand uncertainty exacerbates this centralized tendency and drives model results with higher network connectivity when comparing the SDN to BDN. However, as we can see from the network configurations of SDiN and ISN, accounting for the uncertainty of disruptions results in more decentralized capacity allocation, smaller logistics hub designs and more robust network interconnection. This is because irreplaceable high-capacity hubs are more likely to incur greater penalty costs in hub failures. Finally, the ISN balances the tradeoff between centralization from demand uncertainty and decentralization from disruption uncertainty by establishing 11 logistics hubs across the US Southeast region, with an average daily vehicle throughput capacity of 49.6. It uses relatively few hubs to create a highly connected network where each hub is connected to around 4 other hubs, leading to a resilient and stable network.

Table 2: Basic information of the four optimal hub networks

| Optimal Networks for Each of Stress Testing Levels | BDN | SDN | SDiN | ISN |
|---|---|---|---|---|
| Hub Network Throughput Capacity (Vehicles per day) | 348 | 481 | 485 | 546 |
| Number of Active Hubs | 7 | 9 | 13 | 11 |
| Hub Average Throughput Capacity (Vehicles per day) | 49.7 | 53.4 | 37.3 | 49.6 |
| Hub Connectivity (Average degree of hub nodes) | 2.0 | 3.1 | 4.9 | 4.2 |

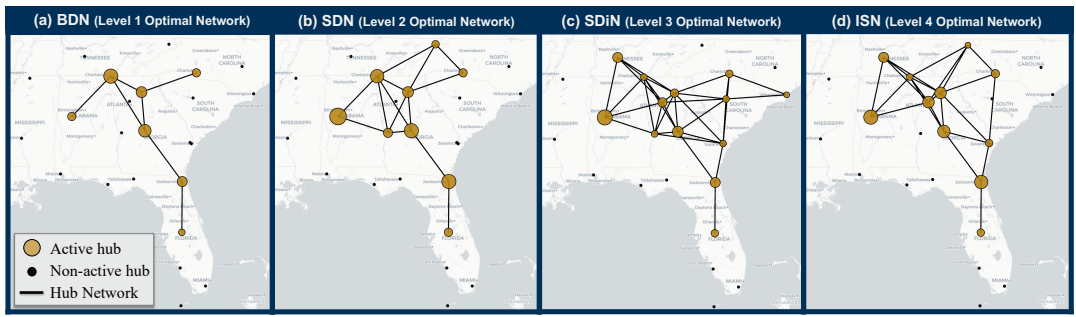

Figure 3: Optimal hub network configuration at each stress testing level

To test the performances of the four hub networks, we perform simulation studies for each optimal network at four stress testing levels. In the simulation, the demand scenarios are the actual daily demands from the dataset, while the disruption scenarios are randomly generated daily. We consider the case instance where the simulation runs for three months, the delivery deadline is one day, and all disruptions last for one day. The simulation's key performance indicators include delivery timeliness and cost-effectiveness. Overall, increasing scenario stress testing levels from 1 to 4 results in decreased on-time delivery rates and increased daily total costs for each network as illustrated in Fig. 4.

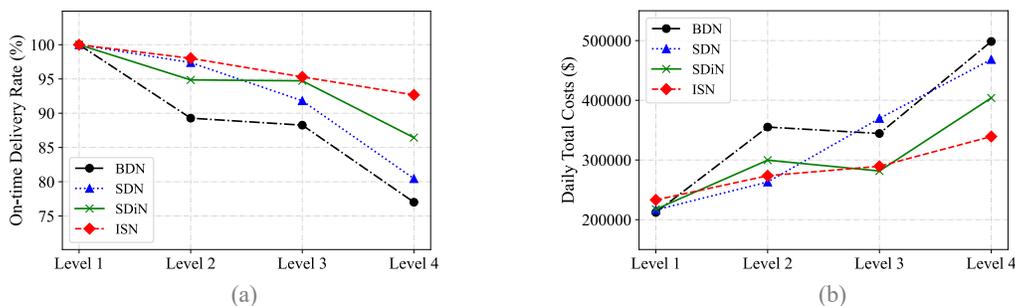

Figure 4: On-time delivery rate and average daily total cost for each network at each stress testing level

Fig. 4 (a) shows that the ISN consistently has the highest on-time delivery rate compared to the other networks and is more resilient to uncertainties since it has a smaller decreasing slope as uncertainty increases. In Fig. 4(b), the optimal hub network at each uncertain level tends to have the lowest cost in its corresponding level, since the goal of our optimization model is to minimize total cost. However, moving from level 1 to level 4 as uncertainty becomes more complex, the ISN is more resilient and robust to keep daily total costs relatively low under $350,000 through all stress testing levels. Other networks, specifically the BDN from the deterministic optimization model, are more sensitive to uncertainties, and their daily costs increase rapidly as uncertainty becomes complex. Overall, our two-stage stochastic optimization model yields different hub network configurations based on various sets of input scenarios. Incorporating more uncertainties, such as demand fluctuations and geographical disruptions, leads to a more resilient network. The aforementioned illustrates the significant impact of capacity deployment on network resilience.

## 6. Conclusion

In this paper, we have proposed a two-stage stochastic optimization model to address spatial-temporal uncertainty in hyperconnected transportation networks, which can be solved efficiently utilizing the SAA method. We further analyzed the separate and combined effects of demand uncertainty and geographical disruptions on hub network configurations by solving the model under different stress testing levels. Our stress testing simulation results demonstrate that our approach provides practical , effective, and resilient hub network capacity deployment solution with high on-time delivery rates, and low total cost under multiple demand and disruption uncertainty scenarios.

Overall, we have proposed an adaptive logistics hub capacity deployment approach in hyperconnected transportation settings. Additionally, this model can be extended to multi-day shipping by modifying each scenario in the input scenario set from a one-day scenario to a k-day scenario and applied to a variety of hub networks that require freight transfer via hubs. With prior knowledge of uncertainties, our model can provide a hub network that effectively adapts to uncertainty and deliver high performances in terms of network resilience, cost savings, and delivery timeliness.

We here propose several research avenues for future research. First, advanced computational methods like Bender's decomposition can be studied to speed up the solution approach and make it more efficient for large-scale instances. Second, beyond our current static hub capacity deployment, it is worthwhile to explore the dynamic capacity deployment by accounting for the changes of throughput flow across time. Third, hub capacity can be further divided into multiple categories, including the loading zone capacity, parking zone capacity, driver-switching zone capacity and so on. Such capacity categorization is usually studied together with hub facility and layout design.